\documentclass[11pt]{article}
\usepackage{amsfonts,latexsym,rawfonts,amsmath,amssymb,amsthm}
\numberwithin{equation}{section}
\newtheorem{theorem}{Theorem}[section]
\newtheorem{lem}[theorem]{Lemma}
\newtheorem{thm}[theorem]{Theorem}

\newtheorem{cor}[theorem]{Corollary}
\def\s{\,\,\,\,}
\def\lan{\langle}
\def\ran{\rangle}
\def\dint{\displaystyle{\int}}
\def\mv{1.7ex}
\title{The $Q$-curvature on a 4-dimensional Riemannian manifold $(M,g)$ with
$\int_MQdV_g=8\pi^2$}
\author{Jiayu Li, Yuxiang Li, Pan Liu}
\date{}
\begin{document}
\maketitle

\section{Introduction}
One of the most important problem in conformal geometry is the
construction of conformal metrics for which a certain curvature
quantity equals a prescribed function, e.g. a constant. In two
dimensions, the problem of prescribed Gaussian curvature asks the
following: given a smooth function $K$ on $(M, g_0)$, can we find
a metric $g$ conformal to $g_0$ such that $K$ is the Gaussian
curvature of the new metric $g$? If let $g = e^{2u} g_0$ for some
$u \in C^{\infty}(M)$, then the problem is equivalent to solving
the nonlinear elliptic equation:
\begin{equation}\label{0010}
\Delta u + K e^{2u} - K_0 = 0,
\end{equation}
where $ \Delta$ denotes the Beltrami-Laplacian of $(M, g_0)$ and
$K_0$ is the Gaussian curvature of $g_0$.

In dimension four, there is an analogous formulation of equation
(\ref{0010}). Let $(M, g)$ be a compact Riemannian four manifold,
and let $Ric$ and $R$ denote respectively the Ricci tensor and the
scalar curvature of $g$. A natural conformal invariant in
dimension four is
\[
Q = Q_g = - \frac{1}{12}(\Delta R - R^2 + 3 |Ric|^2).
\]
Note that, under a conformal change of the metric
\[
\tilde{g} = e^{2u} g,
\]
the quantity $Q$ transforms according to
\begin{equation} \label{002}
2 Q_{\tilde{g}} = e^{- 4u}(P u + 2 Q_g),
\end{equation}
where $P = P_g$ denotes the Paneitz operator with respect to $g$,
introduced in [P]. For any $g$ the operator $P_g$ acts on a smooth
function $u$ on $M$ via
\[
P_g (u) = \Delta_{g}^2 u + div(\frac{2}{3}R_{g}- 2 Ric_{g})du,
\]
which plays a similar role as the Laplace operator in dimension
two. Note that the Paneitz operator is conformal invariant in the
sense that
\[
P_{\tilde{g}} =  e^{-4 u}P_g
\]
for any conformal metric $\tilde{g} = e^{2u} g$.

It follows that the expression  $k = k_g := \int_M Q d V_g$ is
conformally invariant. Moreover, in view of relation (\ref{002}),
a natural problem to propose is to prescribe the $Q$-curvature:
that is, to ask whether on a given four-manifold $(M, g)$ there
exists a conformal metric $\tilde{g}: = e^{2u} g$ for which the
$Q$-curvature of $\tilde{g}$ equlas the prescribed function
$\tilde{Q}$? This is related to solving the following equation
\begin{equation} \label{003}
P_g u + 2 Q_g = 2 \tilde{Q} e^{4u}.
\end{equation}
This equation is the Euler-Language equation of the functional
\begin{equation}
II_{g}(u)=\dint_M  u P_g u dV_g+ 4 \dint_M Q_g u dV_g- (\int_M Q_g
dV_g) \log \dint_M \tilde{Q} e^{4u}dV_g.
\end{equation}

A partial affirmative answer to the problem (\ref{003}) in the
case that $\tilde{Q}$ equals $\mbox{some constant}$ is given by
Chang-Yang \cite{C-Y} provided that the Paneitz operator is weakly
positive and the integral $k$ is less than $8 \pi^2$. In view of a
result of Gursky \cite{G} the former hypothesis is satisfied
whenever $k
> 0$ and provided $(M, g)$ is of positive Yamabe type. The result
of Chang-Yang has been extended recently by Djadli-Malchiodi
\cite{D-M} to the case in which $P_g$ has no kernel and $k$ is not
positive integer multiple of $8 \pi^2$.\\

In the critical case, when $k = 8\pi^2$, the study of equation
(\ref{003}) becomes rather delicate. In this case the functional
$II_{g}$ fails to satisfy standard compactness conditions like the
Palais-Smale condition, and generally blow-up may occur. Note that
when $(M, g)  = (S^4, g_c)$, the above equation (\ref{003}) is
reduced to the following one
\begin{equation} \label{005}
P_g u +  6 = 2 \tilde{Q} e^{4u}.
\end{equation}
This is the analogue of the well-known Nirenberg's problem. This
problem has been recently studied by many authors (please see
\cite{W-X}, \cite{M-St} and the reference there in). We remark
that, similar to Nirenberg's problem, there are some obstructions
to the existence of solution to equation (\ref{005}) in the
standard four-sphere case. The Gauss-Bonnet-Chern formula implies
that there could not be a solution if $\tilde{Q} \leq 0$. On the
other hand, one has the
identities of Kazdan-Warner type to this equation.\\

The main goal of this paper is to study the equation (\ref{003})
with critical value $k = 8 \pi^2$. We shall pursue a variational
approach which was used in \cite{D-J-L-W}.
Let $(M,g)$ be any closed four dimensional Riemannian manifold
with positive $P_g$, i.e., $\int_M uP_gu dV_g\geq 0$ and $ker P_g=\{constants\}$.
Then we have
$$\dint_{M} uP_gu dV_g\geq \lambda\dint_M|\nabla_gu|^2dV_g$$
for some positive $\lambda$ and the following
improved Adams-Fontana inequality \cite{C-Y}:
\begin{equation}\label{fo}
\log\dint_Me^{4u}dV_g \leq \frac{1}{8\pi^2}
\dint_Mu P_gu dV_g + 4 \dint_MudV_g  + C, ~ \forall u \in
W^{2,2}(M).
\end{equation} \\
 We consider (for any small $\epsilon>0$)
$$II_\epsilon(u)=\dint_M \lan u,u\ran dV_g + 4(1-\frac{\epsilon}{8\pi^2})
\dint_MQ_gudV_g-(8\pi^2-\epsilon)\log\dint_M\tilde{Q}
e^{4u}dV_g,$$ where we denote
\[
\lan u,v\ran =\Delta_g u \Delta_g v+(\frac{2}{3}R_g (\nabla
u,\nabla v)- 2Ric_g (\nabla u,\nabla v)).
\]
By using the inequality (\ref{fo}), it is not so difficult to
prove that
$$\inf II_\epsilon(u)>-\infty, \forall \epsilon > 0, ~ \mbox{and moreover,}
~  II_\epsilon ~ \mbox{has a minimal point} ~  u_\epsilon. $$

For this minimizing sequence $u_{\epsilon}$, two possibilities may
occur: let $m_\epsilon=u_\epsilon(x_\epsilon)=\max_{x\in M}
u_\epsilon(x)$,

(1) $\sup\limits_\epsilon m_\epsilon < +\infty$, then, by passing
to a subsequence, $\{u_\epsilon\}$ converges to some $u_0$ as
$\epsilon\rightarrow 0$, and $u_0$ minimizes $II$.

(2) $m_\epsilon\rightarrow +\infty$, as $\epsilon \rightarrow 0$.
We call, in this case, the $u_\epsilon$ blows up.\\

One of the main concern is to prove that, if the second case
happens, then we find an explicit bound for the $II_{\epsilon}$.
More precisely, we have
\begin{equation}\label{007}
\inf_{u\in W^{2,2}(M)}II(u)\geq\Lambda_g(\tilde{Q},p),
\end{equation}
where
$$\Lambda_g(\tilde{Q},p)=-16\pi^2\log\frac{\sqrt{3\tilde{Q}(p)}}{12}-8\pi^2\log8\pi^2
-16\pi^2S_0(p)+2\dint_MQG_pdV_g+(8/3-16)\pi^2,$$ $p$ is the bubble
point, and $S_0(p)$ is the constant term of the
Green function at point $p$ (please see section 6).\\

On the other hand , if we can construct some test function
sequence $\phi_\epsilon$, s.t.
$$II(\phi_\epsilon)<\Lambda_g(\tilde{Q},p),$$
we see that the blow-up does not happen. Therefore, we can get
some sufficient condition under which (\ref{003}) has a
solution.\\

One of our main theorem in this paper is as follows.

\begin{thm}\label{main1}
Let $(M,g)$ be a closed Riemannian manifold of dimension four, with
$k = 8\pi^2$.  Suppose  ${P_g}$ is positive. If the
$\inf\limits_{u\in W^{2,2}(M)}II(u)$ can not be attained, i.e.
equation (\ref{003}) has no minimal solution, then
\begin{equation}\label{007}
\inf_{u\in W^{2,2}(M)}II(u)=\inf_{p\in M}\Lambda_g(\tilde{Q},p).
\end{equation}\\
\end{thm}

Now let $p'$ be a point s.t.
$$\Lambda_g(\tilde{Q},p')=\inf_{x\in M}\Lambda_g(\tilde{Q},x),$$
we will prove that $p'$ is in fact determined by the conformal
class $[g]$ of $(M,g)$.\\

Another main result in this paper is the existence theorem of the
equation (\ref{003}).

\begin{thm}\label{main2}
Let $(M,g)$ be a closed Riemannian manifold of dimension four,
with $k = 8\pi^2$.  Suppose ${P_g}$ is positive. Let $\tilde{Q}$ be a positive smooth
function on $M$. Assume that $\Lambda_g(\tilde{Q},x)$ achieves its
minimum at the point $p'$.  If
$$\tilde{Q}(p')(\Delta_{g}S(p')+4|\nabla_{g}S(p')|^2-\frac{R(p')}{18})+
     [(2\nabla_{g}S\nabla_{g} \tilde{Q})(p')+\frac{1}{4}\Delta_{g}\tilde{Q}(p')]>0,$$
then equation (\ref{003}) has a minimal
solution.\\
\end{thm}

\begin{cor}\label{main3}
With the assumption as in Theorem \ref{main2}. If
$$\Delta_{g}S(p')+4|\nabla_{g}S(p')|^2-\frac{R(p')}{18}
> 0, $$
then $M$ has a constant $Q$-curvature up to conformal
transformations.
\end{cor}

It is interesting to note that, in four-dimensional case, the
method in \cite{D-J-L-W} can not be directly used. In our case
there are some interesting points happens, one is that we use the
method \cite {M-2} to collect the nice information around the
bubble points. The second one is a new technique used in the
derivation of (\ref{007}), where the key point is to calculate
\begin{equation}\label{008}
\dint_{B_\delta\setminus
B_{Lr_\epsilon}(x_\epsilon)}|\Delta_gu_\epsilon|^2dV_g.
\end{equation}
Since the equation (\ref{003}) does not satisfy the Maximal
Principle, the method used in \cite{D-J-L-W} does not work here.
We will apply the capacity to get the lower bound of (\ref{008}).
The usefulness of capacity in similar problems was first
discovered by the second author, and has been used in \cite{Li}
and  \cite{Li-Li}.

We remark that the methods in this paper also work for the
equation
\begin{equation}
\label{012} P_gu + 16\pi^2 = 2he^{4u},
\end{equation}
on any 4-dimensional manifold under the assumptions that $
P_g$ is positive and $Vol=1$. Therefore Theorem \ref{main1} and
Theorem \ref{main2} hold for equation (\ref{012}) (just change
$\tilde{Q}$ to $h$).

\section{Preliminary estimate}

In this section we collect some useful preliminary facts and then
drive some estimates for the solutions. We start with the
following lemma.

\begin{lem}
For any $\epsilon>0$, $II_\epsilon$ has a minimal point.
\end{lem}
\proof By using the inequality (\ref{fo}), it is easy to see that, when $\int_MudV_g=0$, we have
$$\begin{array}{lll}
    II_\epsilon(u)&=&\dint_M u P_gu dV_g+4(1-\frac{\epsilon}{8\pi^2})\dint_MQ u dV_g
                     -(8\pi^2-\epsilon)\log\dint_M\tilde{Q} e^{4u}dV_g\\[\mv]
                  &\geq&C+\frac{\epsilon}{8\pi^2}\dint_Mu P_gu dV_g
                  +4(1-\frac{\epsilon}{8\pi^2})\dint_MQ u dV_g\\[\mv]
                  &\geq&C+\lambda\frac{\epsilon}{8\pi^2}\dint_M|\nabla_gu|^2dV_g
                   +4(1-\frac{\epsilon}{8\pi^2})\dint_MQ u dV_g.
\end{array}$$
For any $\epsilon_1>0$, we have
$$\dint_MQudV_g\leq \epsilon_1\dint_M|u|^2+C_\epsilon\leq \lambda_0\epsilon_1\dint_M
|\nabla u|^2dV_g+C_\epsilon,$$
where $\lambda_0$ is the first eigenvalue of $\Delta$.
Then,
\begin{equation}{\label{n1}}
\dint_M|\nabla_g u|^2dV_g\leq C(\epsilon)II_\epsilon(u)+C
\end{equation}
and then
\begin{equation}{\label{n2}}
\dint_M|\Delta_gu|^2dV_g\leq \frac{8\pi}{\epsilon}II_\epsilon(u)+C.
\end{equation}

Let $u_k = u_{\epsilon, k}$ be a minimizing sequence of
$II_\epsilon$, i.e.
$$II_\epsilon(u_k)\rightarrow\inf II_\epsilon(u)= A,$$
which, together with the above inequality, implies that
\[
\dint_M|\Delta_g u_k|^2dV_g\leq C,
\]
for some constant $C$ which may depend on $\epsilon$. Therefore,
by passing to a subsequence, we have $u_k\rightharpoondown
u_{\epsilon}$ and
$$\dint_M|\Delta_gu_k|^2dV_g\rightarrow B.$$
Since the functional $II_\epsilon$ is invariant under a
translation by a constant,  we may assume that $\int_Mu_kdV_g=0$,
then by (\ref{fo}), we can see that $e^{4u_k}\in L^p$ for any
$p>0$.

Set
$$II_\epsilon(u_k) :=\dint_M|\Delta_g u_k|^2dV_g+\dint_MF(u_k)dV_g,$$
then we have,
$$\lim_{k\rightarrow+\infty}\dint_MF(u_k)dV_g= A-B,\s
and\s \lim_{k\rightarrow+\infty,m\rightarrow+\infty}
\dint_MF(\frac{u_k+u_m}{2})dV_g= A-B.$$ Since
$II_\epsilon(\frac{u_k+u_m}{2})\geq A$, we have
$$\frac{1}{4}\dint_M(|\Delta_gu_k|^2+|\Delta_gu_m|^2)dV_g
+\frac{1}{2}\dint_M\Delta_gu_k\Delta_gu_mdV_g\geq B.$$ Hence
$$\lim_{k\rightarrow+\infty,m\rightarrow+\infty}
\dint_M\Delta_gu_k\Delta_gu_mdV_g\geq B.$$ Then
$$\begin{array}{l}
   \lim\limits_{k\rightarrow+\infty,m\rightarrow+\infty}
      \dint_M|\Delta_g(u_k-u_m)|^2dV_g=\\[\mv]
   \s\s\s\lim\limits_{k\rightarrow+\infty,m\rightarrow+\infty}
        (\dint_M|\Delta_g u_k|^2dV_g+\dint_M|\Delta_gu_m|^2dV_g-2
        \dint_M\Delta_gu_k\Delta_gu_mdV_g)\leq 0.
   \end{array}$$
Therefore, $\{u_k\}$ is a Cauchy sequence in $W^{2,2}(M)$.

$\hfill\Box$\\

\begin{lem}\label{}
We have $\inf II_\epsilon$ is decreasing in $\epsilon$. Moreover,
$$\lim_{\epsilon\rightarrow 0}II_\epsilon=\inf II.$$
\end{lem}

\proof Since $II_{\epsilon}(u+c)=II_{\epsilon}(u)$, we can assume
that $\int_MQ_gudV_g=0$. Therefore
$$II_{\epsilon'}(u)=II_{\epsilon}(u)+(\epsilon-\epsilon')\dint_M\tilde{Q} e^{4u}.$$
Hence, $\inf II_\epsilon$ is decreasing in $\epsilon$, and $\inf
II\leq \inf II_\epsilon$.

Let $\epsilon'=0$, and $II(u_\epsilon)=\inf II_{\epsilon}(u)$. We
have
$$II(u)\geq II_\epsilon(u_\epsilon)-\epsilon\dint_M\tilde{Q} e^{4u}dV_g.$$
Letting $\epsilon\rightarrow 0$, we get that $\inf II\geq
\lim\limits_{\epsilon \rightarrow 0}\inf II_\epsilon$.

$\hfill\Box$\\

Now let $u_\epsilon$ be the minimal point of $II_\epsilon$, it is
clear that $u_\epsilon$ satisfies the following equation:
$$\left\{\begin{array}{l}
           P_gu_\epsilon+2(1-\frac{\epsilon}{8\pi^2})Q_g = 2(1-\frac{\epsilon}{8\pi^2})
                     \tilde{Q} e^{4u_\epsilon}\\[\mv]
          \int_M\tilde{Q} e^{4u_\epsilon}dV_g=8\pi^2.
\end{array}\right.
$$
The same proof of Lemma 2.3 in \cite{M-2} yields the following

\begin{lem}\label{key}
There are constants $C_1(q)$, $C_2(q)$, $C_3(q)$ depending only on
$p$ and $M$ such that, for $r$ sufficiently small and for any
$x\in M$ there holds
$$\dint_{B_r(x)}|\nabla^3u_\epsilon|^qdV_g\leq C_1(q)r^{4-3q},\s
\dint_{B_r(x)}|\nabla^2u_\epsilon|^qdV_g\leq C_2(q)r^{4-2q},$$ and
$$\dint_{B_r(x)}|\nabla u_\epsilon|^qdV_g\leq C_3(q)r^{4-q}$$
where,
respectively, $q<\frac{4}{3}$, $q<2$, and $q<4$.\\
\end{lem}

\section{The proof of Theorem \ref{main1}}
 Let $x_\epsilon$ be the maximum point
of $u_\epsilon$. Assume $m_\epsilon= u_\epsilon(x_\epsilon)$,
$r_\epsilon=e^{-m_\epsilon}$, and $x_\epsilon\rightarrow p$. Let
$\{e_i(x)\}$ be an  orthogonal basis of $TM$ near $p$ and $exp_x:
T_xM\rightarrow M$ be the exponential mapping. The smooth mapping
$E: B_\delta(p)\times B_r\rightarrow M$ is defined as follows,
$$E(x,y)=exp_x(y^ie_i(x)),$$
where $B_r$ is a small ball in $\mathbb{R}^n$. Note that
$E(x,\cdot): T_xM\rightarrow M$ are all differential homeomorphism
if $r$ is sufficiently small.

We set
$$g_{ij}(x,y)=\lan
(exp_x)_*\frac{\partial}{\partial
y^i},(exp_x)_*\frac{\partial}{\partial y^j} \ran_{E(x,y)}.$$ It is
well-known that $g=(g_{ij})$ is smooth, and $g(x,y)=I+O(|y|^2)$
for any fixed $x$. That is, we are able to find a constant $K$,
s.t.
$$\|g(x,y)-I\|_{C^0(B_\delta(p)\times B_r)}\leq K|y|^2$$
when $\delta$ and $r$ are sufficiently small. Moreover, for any
$\varphi\in C^\infty (B_\rho(x_k))$ we have
$$\Delta_gu_\epsilon
=\frac{1}{\sqrt{|g|}}\frac{\partial}{\partial
x^k}(\sqrt{|g|}g^{km} \frac{\partial
u_\epsilon(E(x_\epsilon,x))}{\partial x^m}),\s |\nabla
u_\epsilon|^2 =g^{pq}\frac{\partial
u_\epsilon(E(x_\epsilon,x))}{\partial x^p} \frac{\partial
u_\epsilon(E(x_\epsilon,x))}{\partial x^q},$$ and
$$\dint_{B_\delta(x_k)}\varphi dV_g=\dint_{E^{-1}(x_k,y)B_{\delta(x_k)}}
\varphi(E^{-1}(x_k,y))\sqrt{|g|}dy.$$

We define
$$\tilde{u}_\epsilon(x)=u_\epsilon(E(x_\epsilon,x)),$$
and
$$v_\epsilon(x)=\tilde{u}_\epsilon(r_\epsilon x), v_\epsilon'=v_\epsilon-m_\epsilon.$$
Now $v_\epsilon$, $v_\epsilon'$  are  functions defined on
$B_{\frac{r}{2r_\epsilon}} \subset\mathbb{R}^n$.

We have
\begin{equation}\label{3.1}
\Delta_{g_\epsilon}^2v_\epsilon'=r_\epsilon^2O(|\nabla^2v_\epsilon'|)+
r_\epsilon^3O(\nabla
v_\epsilon')+\tilde{Q}_g(E(x_\epsilon,r_\epsilon
x))e^{4v_\epsilon'}.
\end{equation}
It follows from Lemma \ref{key} that,
$$\|\nabla^2v_\epsilon'\|_{L^q(B_L)}\leq C(L,q)
\hbox{ and }\|\nabla v_\epsilon'\|_{L^q(B_L)}\leq C'(L,q) \hbox{
for any }q\in(1,2).$$ Then (\ref{3.1}) implies that
$$\|\Delta_{g_\epsilon}(\Delta_{g_\epsilon}v_\epsilon')\|_{L^q(B_L)}\leq C'(L).$$
Using the standard elliptic
estimate, we get
$$\|\Delta_{g_k} v_\epsilon'\|_{W^{2,q}(B_L)}\leq C_2(L).$$
The Sobolev inequality then yields that,
$$\|\Delta_{g_\epsilon}v_\epsilon'\|_{L^q(B_L)}\leq C_3(q,L)\s for \s any\s q\in (0,4).$$
We therefore have
$$\|v_\epsilon'\|_{W^{2,q}(B_L)}\leq C_4(L).$$
Hence, by using the standard elliptic estimates, we see that
$v_\epsilon'$ converge smoothly to $w$, which satisfies
$$\Delta_0^2w = 2\tilde{Q}(p)e^{4w}.$$
Moreover, it is easy to check that
$$\dint_{B_L}\tilde{Q}(p)e^{4w}dx\leq 8\pi^2$$
for any $L>0$. By the result of \cite{Lin}, we have

a) $w=-\log(1+\frac{\sqrt{3\tilde{Q}(p)}}{12}|x|^2)$, with
$$\tilde{Q}(p)\dint_{\mathbb{R}^4}e^{4w}dV_g=8\pi^2,$$

or

b) $w$ has the following asymptotic behavior:

$$-\Delta w\rightarrow a>0\hbox{ as }|x|\rightarrow+\infty.$$
We claim that b) does not happen. If it does, then we have
$$\lim_{\epsilon\rightarrow+0}\dint_{B_R}-\Delta_gv_\epsilon\sim\frac{\omega_3}{4}aR^4.$$
However, it follows from Lemma \ref{key} that
$$\dint_{B_R}|\Delta_{g_\epsilon}v_\epsilon'|dV_g\leq CR^2.$$
This shows the case b) does not happen.

For simplicity, let $\lambda=\frac{\sqrt{3Q(p)}}{12}$, so that
we have
$$w=-\log(1+\lambda |x|^2).$$

Now, we consider the convergence of $u_\epsilon$ outside the
bubble. By Lemma \ref{key}, $u_\epsilon$ is bounded in $W^{3,q}$
for any $q<\frac{4}{3}$. Then, it is easy to check that
$u_\epsilon-\bar{u}_\epsilon\rightharpoondown G_p$, where
$$P_gG_p+2Q_g=16\pi^2\delta_p,\s\dint_MG_pdV_g=0.$$
To prove the strong convergence of $u_\epsilon-\bar{u}_\epsilon$,
we first show the following lemma.

\begin{lem} Given $\Omega\subset\subset M\setminus\{p\}$, there holds
$$\dint_\Omega e^{q(u_\epsilon-\bar{u}_\epsilon)}dV_g<C(\Omega, q)$$
for any $q>0$.
\end{lem}

\proof Let $f_\epsilon=\tilde{Q}_ge^{4u_\epsilon}$. For any $x \in
\Omega$, we have the following representation formula,
$$u_{\epsilon}(x)-\bar{u}_{\epsilon} = - \dint_MG(x,y)Q_gdV_{g,y}+\dint_MG(x,y)f_\epsilon.$$

Hence, if let $\Omega_\epsilon=M\setminus
B_{L\epsilon}(x_\epsilon)$, and
$\mu_\epsilon=1/\int_{\Omega_\epsilon}|f|dV_g$, we have, for any
$q'>0$,
$$e^{ q' \mu_\epsilon(u_\epsilon-\bar{u}_\epsilon+\int_M G(x,y)Q_gdV_g)}
=e^{\int_{\Omega_\epsilon} q' G(x,y)\mu_\epsilon
f_\epsilon(y)dV_{g,y}+ \int_{B_{Lr_\epsilon}} q' G(x,y)\mu_\epsilon
f_\epsilon(y)dV_{g,y}}.$$
Notice that for any $x\in\Omega$, we
have
$$\dint_{B_{Lr_\epsilon}(x_\epsilon)}q'|G(x,y)|\mu_\epsilon f_\epsilon(y)dV_{g,y}
\leq C_1(L)\dint_{B_{Lr_\epsilon}(x_\epsilon)}f_\epsilon(y)dV_g
\leq C_2(L),$$
and
$$e^{\int_{\Omega_\epsilon}q' G(x,y)\mu_\epsilon f_\epsilon(y)dV_{g,y}}
\leq\dint_{\Omega_\epsilon}\frac{f_\epsilon(y)}{\|f_\epsilon\|_{L^1(\Omega_\epsilon)}}e^{
q'G(x,y)} dV_{g,y}.$$ Therefore, by using the Jensen's inequality
and the Fubini's theorem, we obtain
$$\begin{array}{lll}
   \dint_\Omega e^{\int_{\Omega_\epsilon}q' G(x,y)\mu_\epsilon f_\epsilon(y)dV_{g,y}}dV_g
        &\leq&\dint_\Omega\frac{f_\epsilon(y)}{\|f_\epsilon\|_{L^1(\Omega_\epsilon)}}
        (\dint_{
          \Omega_\epsilon}e^{q' G(x,y)}
          dV_{g,x})dV_{g,y}\\
        &\leq&C\dint_\Omega\frac{f_\epsilon(y)}{\|f_\epsilon\|_{L^1(\Omega_\epsilon)}}
        (\dint_{
           \Omega_\epsilon}\frac{1}{|x-y|^\frac{q'}{8\pi^2}}
       dV_{g,x})dV_{g,y}.
  \end{array}$$
The last integral is finite provided $q'< 32\pi^2$. Hence, for any
$q > 0$, if $\epsilon$ is sufficiently small so that $q \leq q'
\mu_{\epsilon}$ we have
$$\begin{array}{lll}
   \dint_\Omega e^{q (u_{\epsilon}(x)-\bar{u}_{\epsilon})} dx
   \leq \dint_\Omega e^{q'
   \mu_{\epsilon}(u_{\epsilon}(x)-\bar{u}_{\epsilon})}
   dx \leq C  \dint_\Omega e^{\int_{\Omega_\epsilon} q' G(x,y)\mu_\epsilon
   f_\epsilon(y)dV_{g,y}}dV_g
  \leq C.
  \end{array}$$

$\hfill\Box$\\

As a consequence of the above lemma, we have

\begin{lem}
Let $\Omega\subset\subset M\setminus\{x_0\}$. Then
$u_\epsilon-\bar{u}_\epsilon$ converges to $G_{x_0}$ in
$C^k(\Omega)$ as $\epsilon\to 0$.
\end{lem}

\proof It is easy to see that $\bar{u}_\epsilon<C.$ Then the lemma
follows.

$\hfill\Box$\\

Remark: In $B_{\delta_0}$, we set $p=y_\epsilon$ for any
$\epsilon$. Clearly, $y_\epsilon\rightarrow 0$. Then we also have
$u_\epsilon(E(p,x))-\bar{u}_\epsilon\rightarrow G_{p}(E(p,x))$.
Moreover, we may write
$$G(E(p,x))=-2\log{|x|}+ S_0(p) + S_1(x),$$
where $S_0(p)$ is a constant and $S_1=O(r^{2+\alpha})$. It is easy
to check $\tilde{u}_\epsilon-\bar{u}_\epsilon\rightarrow
G(E(p,x))$ smoothly in $B_{\delta_0}\setminus B_\delta$ for any
fixed
$\delta$.\\

Now, we  estimate the lower bound of $\lim\limits_{\epsilon\rightarrow 0}
\int_M\lan u_\epsilon, u_\epsilon\ran dV_g$. We write
$$
\dint_M\lan u_\epsilon,u_\epsilon\ran dV_g = I_1 + I_2 + I_3,
$$
where $I_1, I_2, I_3$ denote the integrals on $M\setminus
B_\delta(x_\epsilon)$, $B_{Lr_\epsilon}(x_\epsilon)$ and
$B_\delta\setminus B_{Lr_\epsilon}(x_\epsilon)$ (any fixed $L$ and
$\delta$) respectively. We remark that the integral $I_1$, $I_2$
can be easily treated due to the above lemmas. On the other hand,
by Lemma \ref{key}, we have
$$\dint_{B_\delta\setminus B_{Lr_\epsilon}(x_\epsilon)}|\nabla_g u_\epsilon|^2 dV_g
\rightarrow\dint_{B_\delta(p)}|\nabla_gG|^2=O(\delta^2).$$
So, the key point is to calculate
$$\dint_{B_\delta(x_\epsilon)\setminus B_{Lr_\epsilon}(x_\epsilon)}
|\Delta_g u_\epsilon|^2dV_g.$$

We are going to prove the following lemma.
\begin{lem}\label{l3.3} We have
$$\dint_{B_\delta(x_\epsilon)\setminus B_{Lr_\epsilon}(x_\epsilon)}
|\Delta_g u_\epsilon|^2dV_g \geq\dint_{B_\delta\setminus
B_{Lr_\epsilon}}|(1-B|x|^2)
\Delta_0\tilde{u}_\epsilon|^2dx+J(L,\epsilon,\delta),$$ for some
$B>0$, where
$$\lim_{\delta\rightarrow 0}\lim_{\epsilon\rightarrow 0}J(L,\epsilon,\delta)=0.$$
\end{lem}
\proof Since we have
$$\begin{array}{lll}
   |\Delta_gu_\epsilon|^2&=&|g^{km}\frac{\partial^2\tilde{u}_\epsilon}{\partial x^k
   \partial x^m}
       +O(|\nabla \tilde{u}_\epsilon|^2)|^2\\[\mv]
   &=&|g^{km}\frac{\partial^2\tilde{u}_\epsilon}{\partial x^k\partial x^m}|^2+
      O(|\nabla^2\tilde{u}_\epsilon|(|\nabla \tilde{u}_\epsilon|)
       )+O((|\nabla \tilde{u}_\epsilon|^2)),
\end{array}$$
and since $\tilde{u}_\epsilon-\bar{u}_\epsilon$ converges to
$G_{p}(E(p,x))$ in $W^{3,q}$ for any $q<\frac{4}{3}$, we get
$$\begin{array}{l}
   \int_{B_\delta\setminus B_{Lr_\epsilon}}O(|\nabla^2\tilde{u}_\epsilon|(|\nabla
   \tilde{u}_\epsilon|)+O(|\nabla \tilde{u}_\epsilon|^2)\\[\mv]
     \s\s\s\s\s\s\leq C(
      \|\nabla^2G_{p}\|_{L^q(B_\delta\setminus B_{Lr_\epsilon})}
      \|\nabla_gG_{p}\|_{L^{q'}}
      (B_\delta\setminus B_{Lr_\epsilon})+\|G_{p}\|_{W^{1,2}(B_\delta\setminus
        B_{Lr_\epsilon})})\\[\mv]
     \s\s\s\s\s\s= J(L,\epsilon,\delta),
  \end{array}$$
where $\frac{3}{2}<q<2$, and $\frac{1}{q'}+\frac{1}{q} = 1$ .

Let $g^{km}=\delta^{km}+A^{km}$, with $|A^{km}|\leq K|x|^2$ for
any $\epsilon, k, m$. Consequently we have
$$|g^{km}\frac{\partial^2\tilde{u}_\epsilon}{\partial x^k\partial x^m}|^2=
|\Delta_0\tilde{u}_\epsilon|^2+2\sum_{s,t}A^{st}\Delta_0
\tilde{u}_\epsilon
\frac{\partial^2\tilde{u}_\epsilon}{\partial x^s\partial
x^t}+\sum_{k,m,s,t} A^{km}A^{st}\frac{\partial^2
\tilde{u}_\epsilon} {\partial
 x^k\partial x^m}
\frac{\partial^2\tilde{u}_\epsilon}{\partial x^s\partial x^t}.$$
It is clear that
$$2\dint_{B_\delta\setminus B_{Lr_\epsilon}}|A^{st}\Delta_0
\tilde{u}_\epsilon
\frac{\partial^2\tilde{u}_\epsilon}{\partial x^s\partial x^t}|\leq
K\dint_{B_\delta\setminus B_{Lr_\epsilon}} (|x|^2|\Delta_0
\tilde{u}_\epsilon|^2+|x|^2|
\frac{\partial^2\tilde{u}_\epsilon}{\partial x^s\partial
x^t}|^2)dx,$$
and
$$\begin{array}{lll}
   \dint_{B_\delta\setminus B_{Lr_\epsilon}}|x|^2|
     \frac{\partial^2\tilde{u}_\epsilon}{\partial x^s\partial x^t}|^2dx&=&
     \dint_{B_\delta\setminus B_{Lr_\epsilon}}|x|^2\frac{\partial^2\tilde{u}_\epsilon}
      {\partial x^t\partial x^t}
      \frac{\partial^2\tilde{u}_\epsilon}{\partial x^s\partial x^s}dx
      +\dint_{B_\delta\setminus B_{Lr_\epsilon}}
     O(|x|\,|\nabla \tilde{u}_\epsilon|\,|\nabla^2 \tilde{u}_\epsilon|)dx\\[\mv]
   &&+\dint_{\partial(B_\delta\setminus B_{Lr_\epsilon})}
     |x|^2\frac{\partial \tilde{u}_\epsilon}
      {\partial x^t}
      \frac{\partial^2\tilde{u}_\epsilon}{\partial x^s\partial x^t}\lan\frac{\partial}
      {\partial x^t},
      \frac{\partial}{\partial r}\ran ds\\[\mv]
   &&+\dint_{\partial(B_\delta\setminus B_{Lr_\epsilon})}
       |x|^2\frac{\partial \tilde{u}_\epsilon}
      {\partial x^t}
      \frac{\partial^2\tilde{u}_\epsilon}{\partial x^s\partial x^s}\lan\frac{\partial}
      {\partial x^s},
      \frac{\partial}{\partial r}\ran)ds\\[\mv]
   &=&\dint_{B_\delta\setminus B_{Lr_\epsilon}}|x|^2\frac{\partial^2\tilde{u}_\epsilon}
      {\partial x^t\partial x^t}
      \frac{\partial^2\tilde{u}_\epsilon}{\partial x^s\partial
      x^s}dx+J(L,\epsilon,\delta).
  \end{array}$$
Hence,
$$2\sum_{k,s,t}\dint_{B_\delta\setminus B_{Lr_\epsilon}}|A^{st}\Delta_0
\tilde{u}_\epsilon
\frac{\partial^2\tilde{u}_\epsilon}{\partial x^s\partial x^t}|\leq
4K\dint_{B_\delta\setminus
B_{Lr_\epsilon}}|x|^2|\Delta_0\tilde{u}_\epsilon|^2dx+J(L,\epsilon,\delta).$$
A similar argument as above then gives,
$$\dint_{B_\delta\setminus B_{Lr_\epsilon}}\sum_{k,m,s,t}
A^{km}A^{st}\frac{\partial^2 \tilde{u}_\epsilon} {\partial
 x^k\partial x^m}
\frac{\partial^2\tilde{u}_\epsilon}{\partial x^s\partial x^t} \leq
K^2\dint_{B_\delta\setminus
B_{Lr_\epsilon}}|x|^4|\Delta_0\tilde{u}_\epsilon|^2dx+J(L,\epsilon,\delta).$$
This proves the Lemma.

$\hfill\Box$\\

\begin{lem}\label{l3.4}
There is a function sequence $U_\epsilon\in
W^{2,2}(B_\delta\setminus B_{Lr_\epsilon})$ s.t.
$$U_\epsilon|_{\partial B_\delta}=-2\log\delta+S_0(p)+\bar{u}_\epsilon,\s
U_\epsilon|_{\partial
B_{Lr_\epsilon}}=w(L)+m_\epsilon$$
$$\frac{\partial U_\epsilon}{\partial r}|_{
\partial B_\delta}=-\frac{2}{\delta},\s
\frac{\partial U_\epsilon}{\partial r}|_{\partial
B_{Lr_\epsilon}}=w'(L)$$
and
$$\dint_{B_\delta\setminus B_{Lr_\epsilon}}|\Delta_0(1-B|x|^2)(U_\epsilon-
\bar{u}_\epsilon)|^2dx
=\dint_{B_\delta\setminus
B_{Lr_\epsilon}}|(1-B|x|^2)\Delta_0\tilde{u}_\epsilon|^2dx
+J(L,\epsilon,\delta).$$
\end{lem}

\proof

Let $u_k'$ be the solution of
$$\left\{\begin{array}{l}
             \Delta_0^2u_\epsilon'=\Delta_0^2v_\epsilon\\[\mv]
             \frac{\partial u_\epsilon'}{\partial n}|_{\partial B_{2L}}=
                \frac{\partial v_\epsilon}{\partial n}|_{\partial B_{2L}},\s
                 u_\epsilon'|_{\partial B_{2L}}=v_\epsilon|_{\partial B_{2L}}\\[\mv]
             \frac{\partial u_\epsilon'}{\partial n}|_{\partial B_{L}}=
                \frac{\partial w}{\partial n}|_{\partial B_{L}},\s
                 u_\epsilon'|_{\partial B_{L}}=m_\epsilon+w|_{\partial
                 B_{L}}.
          \end{array}\right.$$
We set
$$U_\epsilon'=\left\{\begin{array}{ll}
                        u_\epsilon'(\frac{x}{r_\epsilon})&Lr_\epsilon\leq
                          |x|
                          \leq 2Lr_\epsilon\\[\mv]
                        \tilde{u}_\epsilon(x)&2Lr_\epsilon\leq
                          |x|.
                     \end{array}\right.$$
It is easy to see that $u_\epsilon'-m_\epsilon$ converges to $w$ smoothly on
$B_{2L}\setminus B_L$,
we have
$$\lim_{\epsilon\rightarrow 0}\dint_{B_{2Lr_\epsilon}\setminus B_{Lr_\epsilon}}
(1-B|x|^2)^2(|\Delta_0U_\epsilon'|^2-|\Delta_0\tilde{u}_\epsilon|^2)dx=0.$$

Let $\eta$ be a smooth function which satisfies:
$$\eta(t)=\left\{\begin{array}{ll}
                      1&t\leq 1/2\\[\mv]
                      0& t>2/3
                 \end{array}\right.$$
Set $G_\epsilon=\eta(\frac{|x|}{\delta})(\tilde{u}_\epsilon-S_0(p)
+2\log{|x|^2}-\bar{u}_\epsilon) -2\log{|x|^2}+ S_0(p)$. Recall that
$u_\epsilon-\bar{u}_\epsilon$
converges to $G_p$ smoothly on $M\setminus B_{\frac{\delta}{2}}(p)$, we have
$$G_\epsilon\rightarrow
-2\log{|x|^2}+ S_0(p) + \eta(\frac{|x|}{\delta}) S_1(x),\s
\tilde{u}_\epsilon - G_\epsilon -\bar{u}_\epsilon\rightarrow
(\eta(\frac{|x|}{\delta})-1)S_1(x).$$  Therefore
$$\begin{array}{l}
   \lim\limits_{\epsilon\rightarrow 0}
    \left|\dint_{B_\delta\setminus B_{\delta/2}}|\Delta_0\tilde{u}_\epsilon|^2dx-
    \dint_{B_\delta\setminus B_{\delta/2}}|\Delta_0G_\epsilon|^2dx\right|\\[\mv]
    \s\s\s \leq\sqrt{\int_{B_\delta\setminus B_{\delta/2}}|\Delta_0(\eta(\frac{|x|}
       {\delta})-1)S_1(x)|^2dx
       \int_{B_\delta\setminus B_{\delta/2}}|\Delta_0(G_{p}-2
     \log{|x|^2}+\eta(\frac{|x|}{\delta})S_1(x))|^2dx}\\[\mv]
    \s\s\s\leq C\sqrt{|\log\delta|}\sqrt{\int_{B_\delta\setminus B_{\delta/2}}|\Delta_0
     \eta(\frac{|x|}{\delta})
    S_1(x)|^2dx}\\[\mv]
    \s\s\s\leq C\sqrt{\delta|\log\delta|}.
\end{array}$$
Now set
$$U_\epsilon=\left\{\begin{array}{ll}
                        U_\epsilon'(x)&|x|\leq\frac{\delta}{2}\\[\mv]
                        G_\epsilon(x)+\bar{u}_\epsilon&\delta/2\leq
                          |x|
                          \leq \delta.
                     \end{array}\right.$$
We then have,
$$\begin{array}{lll}
  \dint_{B_\delta\setminus B_{L\epsilon}}|(1-B|x|^2)
   \Delta_0(U_\epsilon-\bar{u}_\epsilon)|^2dx&=&\dint_{B_\delta\setminus B_{Lr_\epsilon}}
    |\Delta_0(1-B|x|^2)(U_\epsilon-\bar{u}_\epsilon)|^2dx\\[\mv]
   &&+ \dint_{B_\delta\setminus
     B_{Lr_\epsilon}} O(|\nabla
        U_\epsilon|^2+|U_\epsilon-\bar{u}_\epsilon|^2)dV_g.
  \end{array}$$

It is easy to check that
$\|U_\epsilon-\bar{u}_\epsilon-G_p(E(p,x))\|_{W^{1,2}
(B_\delta\setminus B_{Lr_\epsilon})}\rightarrow 0$ as
$\epsilon\rightarrow 0$. Therefore, we proved the lemma.

$\hfill\Box$\\

Now, we are going to apply the capacity to derive the lower bound
of
$$\dint_{B_\delta\setminus B_{Lr_\epsilon}}|\Delta_0(1-B|x|^2)(U_\epsilon-\bar{u}_\epsilon)
|^2dx.$$ First we need to calculate
$$\inf_{\Phi|_{\partial B_r}=P_1,\Phi|_{\partial B_R}=P_2,
\frac{\partial\Phi}{\partial r}|_{\partial B_r}=Q_1,
\frac{\partial\Phi}{\partial r}|_{\partial B_R}=Q_2}
\dint_{B_R\setminus B_r}|\Delta_0\Phi|^2dx,$$ where $P_1$, $P_2$,
$Q_1$, $Q_2$ are constants. Obviously, the minimum can be attained
by the function $\Phi$ which satisfies
$$\left\{\begin{array}{l}
             \Delta_0^2\Phi=0\\[\mv]
             \Phi|_{\partial B_r}=P_1\,\,\, ,\Phi|_{\partial B_R}=P_2\,\,\, ,
               \frac{\partial\Phi}{\partial r}|_{\partial B_r}=Q_1\,\,\, ,
                \frac{\partial\Phi}{\partial r}|_{\partial B_R}=Q_2
         \end{array}\right.$$
Clearly, we can set
$$\Phi=A\log{r}+Br^2+\frac{C}{r^2}+D,$$
where $A$, $B$, $C$, $D$ are all constants. Then we have
$$\left\{\begin{array}{l}
              A\log{r}+Br^2+\frac{C}{r^2}+D=P_1\\[\mv]
              A\log{R}+BR^2+\frac{C}{R^2}+D=P_2\\[\mv]
              \frac{A}{r}+2Br-2\frac{C}{r^3}=Q_1\\[\mv]
              \frac{A}{R}+2BR-2\frac{C}{R^3}=Q_2.
          \end{array}\right.$$
We have
$$\left\{\begin{array}{l}
            A=\frac{P_1-P_2+\frac{\varrho}{2} rQ_1+\frac{\varrho}{2} RQ_2}
             {\log{r/R}+\varrho}\\[\mv]
            B=\frac{-2P_1+2P_2-rQ_1(1+\frac{2r^2}{R^2-r^2}\log{r/R})
                +RQ_2(1+\frac{2R^2}{R^2-r^2}\log{r/R})}{4(R^2+r^2)(\log{r/R}+\varrho)},
          \end{array}\right.$$
where $\varrho=\frac{R^2-r^2}{R^2+r^2}$. Furthermore,
$$\dint_{B_R\setminus B_r}|\Delta_0\Phi|^2dx=-8\pi^2 A^2\log{r/R}+
       32\pi^2 AB(R^2-r^2)+32\pi^2B^2(R^4-r^4)$$

In our case, $R=\delta$, $r=Lr_\epsilon$,
$P_1=m_\epsilon-\bar{u}_\epsilon+w(L)+O(r_\epsilon\bar{u}_\epsilon)$,
$P_2= -2\log\delta+S_0(p)+O(\delta\log\delta)$, $Q_1=\frac{2\lambda
L}{r_\epsilon(1+\lambda L^2)}$,
$Q_2=-\frac{2}{\delta}+O(\delta\log\delta)$. If we define
$$\begin{array}{lll}
       N(L,\epsilon,\delta)&=&w(L)+2\log\delta-S_0-\frac{\varrho}{2} \frac{2\lambda L^2}
                              {1+\lambda L^2}\\[\mv]
       &=&w(L)+2\log\delta-S_0-2+O(\delta\log\delta)+O(\frac{1}{L^2})+O(Lr_\epsilon),
  \end{array}$$
and
$$P=\log\delta-\log{L},$$
then we have
$$\begin{array}{lll}
    A^2(-\log{Lr_\epsilon/\delta})&=&(\frac{m_\epsilon-\bar{u}_\epsilon
     +N(L,\epsilon,\delta)
     }{m_\epsilon+P-\varrho})^2(m_\epsilon+P)\\[\mv]
    &=&(1+\frac{P-\varrho}{m_\epsilon})^{-2}(1+\frac{P}
    {m_\epsilon})
     m_\epsilon(1-\frac{\bar{u}_\epsilon}{m_\epsilon}+\frac{N(L,\epsilon,\delta)}
     {m_\epsilon})^2\\[\mv]
    &=&(1-2\frac{P-\varrho}{m_\epsilon}+O(\frac{1}{m_\epsilon^2}))
      (1+\frac{P}{m_\epsilon})m_\epsilon\\[\mv]
    &&\left[(1-\frac{\bar{u}_\epsilon}{m_\epsilon})^2+
        2(1-\frac{\bar{u}_\epsilon}{m_\epsilon})\frac{N(L,\epsilon,\delta)}{m_\epsilon}+
      O(\frac{1}{m_\epsilon^2})+O(e^{-m_\epsilon} m_\epsilon)\frac{\bar{u}_\epsilon}
       {m_\epsilon}\right]\\[\mv]
    &=&m_\epsilon(1-\frac{\bar{u}_\epsilon}{u_\epsilon})^2+2(1-\frac{\bar{u}_\epsilon}
      {m_\epsilon})N(L,\epsilon,\delta)-
    (P-2\varrho)(1-\frac{\bar{u}_\epsilon}{m_\epsilon})^2\\[\mv]
    &&+O(\frac{1}{m_\epsilon})(1-
      \frac{\bar{u}_\epsilon}{m_\epsilon})^2+O(\frac{1}{m_\epsilon}),
  \end{array}$$
and
$$A=-\frac{m_\epsilon-\bar{u}_\epsilon+N(L,\epsilon,\delta)}{m_\epsilon-\log{L}+\log\delta
+\varrho}
=-(1-O(\frac{1}{m_\epsilon}))^{-1}(1-\frac{\bar{u}_\epsilon}{m_\epsilon}
+O(\frac{1}{m_\epsilon})
)=-1+\frac{\bar{u}_\epsilon}{m_\epsilon}+O(\frac{1}{m_\epsilon}).$$
Notice that  $r_\epsilon m_\epsilon\rightarrow 0$ as
$\epsilon\rightarrow 0$, we have
$$\begin{array}{lll}
    B&=&\frac{-2m_\epsilon+2\bar{u}_\epsilon+O(1)+(2\frac{2\delta^2}{\delta^2-
    (Lr_\epsilon)^2
    }+O(\delta\log\delta))
      m_\epsilon}{4(\delta^2+(Lr_\epsilon)^2)(\log{L}-m_\epsilon-\log\delta+\varrho)}\\[\mv]
    &=&-\frac{1}{2\delta^2}(1+\frac{\bar{u}_\epsilon}{m_\epsilon}+O(\frac{1}{
    m_\epsilon}))
      (1-O(\frac{1}{m_\epsilon}))^{-1}\\[\mv]
    &=&-\frac{1}{2\delta^2}(1+\frac{\bar{u}_\epsilon}{m_\epsilon}+O(\frac{1}{
    m_\epsilon})).
  \end{array}$$
It concludes that
$$\begin{array}{lll}
   \dint_{B_\delta\setminus B_{Lr_\epsilon}}|\Delta_0(1-B|x|^2)(U_\epsilon-
   \bar{u}_\epsilon)|^2dx
     &\geq&8\pi^2 m_\epsilon(1-\frac{\bar{u}_\epsilon}{m_\epsilon})^2+16\pi^2
      (1-\frac{\bar{u}_\epsilon}{m_\epsilon})N(L,\epsilon,
     \delta)\\[\mv]
     &&-
     8\pi^2(P-2\varrho)(1-\frac{\bar{u}_\epsilon}{m_\epsilon})^2\\[\mv]
     &&+16\pi^2(1-\frac{\bar{u}_\epsilon}{m_\epsilon})(1+\frac{\bar{u}_\epsilon}{m_\epsilon})
         +8\pi^2(1+\frac{\bar{u}_\epsilon}{m_\epsilon})^2\\[\mv]
     &&+O(\frac{1}{m_\epsilon})(1-
      \frac{\bar{u}_\epsilon}{m_\epsilon})^2
     +O(\frac{1}{m_\epsilon})+J_6(L,\epsilon,\delta).
   \end{array}$$
Using the fact that $\bar{u}_\epsilon\leq C$, we have
$$(8\pi^2-\epsilon)\bar{u}_\epsilon>8\pi^2\bar{u}_\epsilon+\epsilon C.$$
Therefore
$$\begin{array}{lll}
   II_\epsilon(u_\epsilon)
   &\geq&\int_{B_{Lr_\epsilon}(x_\epsilon)}|\Delta_gu_\epsilon|^2dV_g+
    \int_{B_\delta\setminus B_{Lr_\epsilon}}|\Delta_0(1-|B|^2)(U_\epsilon-
    \bar{u}_\epsilon)|^2dx+8\pi^2
    \bar{u}_\epsilon\\[\mv]
   &&+\int_{M\setminus B_\delta(x_0)}\lan G_p,G_p\ran   +4\int_M\tilde{Q}G_pdV_g+J(L,
   \epsilon,\delta)\\[\mv]
   &\geq& 8\pi^2 (m_\epsilon+C_1)(1+\frac{\bar{u}_\epsilon}{m_\epsilon})^2+C_2(1+\frac{
   \bar{u}_\epsilon}{m_\epsilon})+C_3.
\end{array}$$
where $C_1$, $C_2$, $C_3$ are some constants. Note that since
$II_\epsilon(u_\epsilon) < \infty$, we must have
$(1+\frac{\bar{u}_\epsilon}{ m_\epsilon})\rightarrow 0$ as
$\epsilon \rightarrow 0$, i.e.
$\frac{\bar{u}_\epsilon}{m_\epsilon}\rightarrow -1$.

Consequently we have
\begin{equation}\label{3.2}
\begin{array}{l}
   \int_{B_\delta\setminus B_{Lr_\epsilon}}|\Delta_0(1-B|x|^2)(U_\epsilon-
   \bar{u}_\epsilon)|^2dx+
   8\pi^2\bar{u}_\epsilon\\[\mv]
     \s\s\geq8\pi^2 m_\epsilon(1+\frac{\bar{u}_\epsilon}{m_\epsilon})^2+16\pi^2 N(L,
     \epsilon,\delta)(1-\frac{\bar{u}_\epsilon}{m_\epsilon})-8\pi^2(\log\delta-\log{L}
       -2\varrho)(1-\frac{\bar{u}_\epsilon}{m_\epsilon})^2\\[\mv]
     \s\s\s\s+J(L,\epsilon,\delta)\\[\mv]
   \s\s\geq 16\pi^2(1-\frac{\bar{u}_\epsilon}{m_\epsilon})
     N(L,\epsilon,\delta)-8\pi^2(\log\delta-\log{L}
       -2\varrho)(1-\frac{\bar{u}_\epsilon}{m_\epsilon})^2+J(L,\epsilon,\delta).
   \end{array}
\end{equation}
Since we have
$$\Delta_0w=\frac{4\lambda^2|x|^2}{(1+\lambda |x|^2)^2}-
\frac{8\lambda}{1+\lambda|x|^2},$$
a direct calculation yields
that
\begin{equation}\label{3.3}
\dint_{B_L}|\Delta_0w|^2dx
   =16\pi^2\log(1+\lambda L^2)+\frac{8\pi^2}{3}
     +O(\frac{\log L}{L^2}).$$
On the other hand, it is obvious to see that,
$$\dint_{B_\delta(x_\epsilon)}|\nabla u_\epsilon|^2\rightarrow\dint_{B_\delta(x_\epsilon)}
|\nabla G_p|^2=O(\delta\log\delta),
\end{equation}
and
\begin{equation}\label{3.4}
\begin{array}{lll}
    \dint_{M\setminus B_\delta(x_0)}\lan G_p,G_p\ran dV_g&=&
      \dint_{M\setminus B_\delta(x_0)} G_p P_g G_p dV_g -\dint_{\partial B_\delta}
       \frac{\partial G_p}{\partial r}\Delta_g G_p dV_g+\dint_{\partial B_\delta}
       G_p\frac{\partial \Delta G_p}{\partial r} dV_g\\[\mv]
    &&+\dint_{\partial B_\delta}(\frac{2}{3}RG\frac{\partial G}{\partial
    r}-2GRic(dG,dr))dS_g\\[\mv]
    &=&-2\dint_MQ_gG_pdV_g-16\pi^2+16\pi^2(-2\log\delta+S_0(p))+O(\delta\log\delta).
  \end{array}
\end{equation}

Together with Lemma \ref{l3.3}, Lemma \ref{l3.4}, (\ref{3.2}),
(\ref{3.3}) and (\ref{3.4}), we have
$$\begin{array}{lll}
    \lim\limits_{\epsilon\rightarrow 0}II_\epsilon&\geq&
       32\pi^2\lim\limits_{\epsilon\rightarrow 0}N(L,\epsilon,\delta)-32\pi^2
     (\log\delta-\log{L}-2)+16\pi^2\log(1+\lambda L^2)\\[\mv]
    &&+\frac{8\pi^2}{3}+(-2\log\delta+S_0(p))16\pi^2+2\dint_MQ_gG_pdV_g-8\pi^2\log
    8\pi^2
    +O(\delta\log\delta)+O(\frac{\log L}{L^2})\\[\mv]
    & = & -16 \pi^2 \log \frac{1 + \lambda L^2}{L^2} + \frac{8\pi^2}{3} - 16\pi^2 S_0(p)
       - 16\pi^2+2\dint_MQ_gG_pdV_g-8\pi^2\log
    8 \pi^2 \\[\mv]
   &&+O(\delta\log\delta)+O(\frac{\log L}{L^2}).
  \end{array}$$

Letting first $\delta\rightarrow 0$, then $L\rightarrow+\infty$,
we get
$$\lim_{\epsilon\rightarrow 0}II_\epsilon\geq -16\pi^2\log\lambda-8\pi^2\log8\pi^2
-16\pi^2S_0+(8/3-16)\pi^2+2\dint_MQ_gG_pdV_g.$$ This shows the
first part of Theorem \ref{main1}, that is
$$\inf_{u\in W^{2,2}(M)}II(u)\geq\inf_{p\in M}\Lambda_g(\tilde{Q},p). $$ The second part
$$\inf_{u\in W^{2,2}(M)}II(u)\leq\inf_{p\in M}\Lambda_g(\tilde{Q},p)
$$ follows from the proof of Theorem \ref{main2}
in next section.

To end this section, we will prove a conformal property
of $\Lambda_g(\tilde{Q},p)$.

\begin{lem}\label{l3.5}
Let $\tilde{g} \in [g]$: $\tilde{g} = e^{2v}g$ for some $v \in
C^{\infty}(M)$, we have
$$II_{\tilde{g}}(u)=II_g(u + v)-\dint_M\lan v, v\ran dV_g.$$
If we set
$$P_{\tilde{g}}\tilde{G}_y + 2Q_{\tilde{g}}= 16\pi^2\delta_y,$$
then $\tilde{G}_y = G_y - v$. Moreover, for any $y$, we have
$$2\dint_MQ_{\tilde{g}}\tilde{G}_{y}dV_{\tilde{g}}-16\pi^2\tilde{S}_0(y)
=2\dint_MQ_{g}G_{y}dV_g-16\pi^2S_0(y)-\dint_M\lan v,v\ran dV_g.$$
\end{lem}

\proof Since $P_{\tilde{g}}=e^{-4v}P_g$, $2Q_{\tilde{g}}=e^{-4
v}(P_g v + 2Q_g)$, we get
$$\begin{array}{lll}
  II_{\tilde{g}}(u)&=&\dint_M\lan u,u\ran dV_g+2\dint_M(P_gv+2Q_g)udV_g
     -8\pi^2\log\dint_M\tilde{Q}e^{4(u+v)}dV_g\\[\mv]
  &=&\dint_M\lan u+v,u+v\ran dV_g+4\dint_MQ_gudV_g
   -8\pi^2\log\dint_M\tilde{Q}e^{4(u+v)}dV_g-\dint_M\lan v,v\ran dV_g\\[\mv]
  &=&II_g(u+v)-\dint_M\lan v,v\ran dV_g.
  \end{array}$$
On the other hand, we have
$$P_{\tilde{g}}(G-v)+2Q_{\tilde{g}}=e^{-4v}(P_gG+2Q_g)=16\pi^2e^{-4v}\delta_{y,g}
=16\pi^2\delta_{y,\tilde{g}}.$$ Since
$dist_{\tilde{g}}(y,x)=e^{v(y)}dist_g(y,x)+O(dist_g(y,x))^2$,
we have
$$\begin{array}{lll}
   \tilde{G}_y&=&G_y-v\\[\mv]
     &=&-2\log dist_g(y,x)+S_0(y) - v(y) +O(dist(y,x))\\[\mv]
     &=&-2\log dist_{\tilde{g}}(y,x)+v(y)+S_0(y)+O(dist(y,x)).
   \end{array}$$
Thus $\tilde{S}_0(y)=S_0(y)+v(y)$. Moreover, we have
$$\begin{array}{lll}
  \dint_MQ_{\tilde{g}}\tilde{G}_ydV_{\tilde{g}}&=&
     \dint_M(P_gv+2Q_g)(G_y-v)dV_g\\[\mv]
  &=&(\dint_MG_yP_gvdV_g+2\dint_MQ_gvdV_g)+2\dint_MQ_gG_ydV_g-\dint_Mv P_gv dV_g\\[\mv]
  &=&16\pi^2v(y)+2\dint_MQ_gG_ydV_g-\dint_M v P_g v dV_g,
  \end{array}$$
this proves the lemma.

 $\hfill\Box$\\

\section{Testing function}
In this section we will construct a blow up sequence
$\phi_\epsilon$ s.t.
$$II(\phi_\epsilon)<\inf_{x\in M}\Lambda(x).$$

We use standard notation from [L-P]. In a local coordinate system
$\{x^i\}$, we denote
$$R_{ijkl}= <R(\partial_k,\partial_l)\partial_j, \partial_i>, \s R_{ij}=-g^{jk}R_{ijkl},$$
where $R$ is the curvature operator, defined as follows,
$$R(X,Y)= \nabla_X\nabla_Y -\nabla_Y\nabla_X-\nabla_{[X,Y]}$$

Suppose that $p'$ is a point such that $\Lambda(p')=\inf_{x\in
M}\Lambda(x)$.

We know that, locally we have
$$g_{pq} = \delta_{pq}+\frac{1}{3}R_{pijq}(p')x^ix^j+
\frac{1}{6}R_{pijq,k}(p')x^ix^jx^k+
(\frac{1}{20}R_{pijq,kl}
+\frac{2}{45}R_{pijm}(p')R_{qklm}(p'))x^ix^jx^kx^l+O(r^5).$$
$$|g|=1-\frac{1}{3}R_{ij}x^{ij}-\frac{1}{6}R_{ij,k}(p')x^{ijk}
-(\frac{1}{20}R_{ij,kl}(p')+
\frac{1}{90}R_{hijm}(p')R_{hklm}(p'))x^ix^jx^kx^m+O(r^5)$$

In the sequel, let us denote
$$x^{i_1\cdots i_m}_{j_1\cdots j_n}=x^{i_1\cdots i_m j_1\cdots j_n},\s and\s
\alpha^{i_1\cdots i_m}_{j_1\cdots j_n}=\frac{1}{2\pi^2}\dint_{S^3}
x^{i_1\cdots i_m j_1\cdots j_n}ds,$$ then around the point $p'$ we
write
$$g^{km}=\delta^{km}+M^{km}=\delta^{km}+M_{km}^{ij}x^{km}+M_{kms}^{ij}x^{kms}
+M_{kmst}^{ij}x^{kmst}+O(r^5)$$
$$M=M^{ij}\delta_{ij}=M_{km}x^{km}+M_{kms}x^{kms}+M_{kmst}x^{kmst}+O(r^5),$$
$$ \sqrt{|g|}=1-\frac{1}{6}R_{ij}x^{ij}+K_{ijk}x^{ijk}+K_{ijkm}x^{ijkm}+O(r^5).$$
$$N^k = -g^{ij}\Gamma_{ij}^k
=N^k_{i}x^{i}+N^k_{ij}x^{ij}+N^{k}_{ijm}x^{ijm} + O(r^5).$$
It is easy to check that $M_{km}^{ij}=-\frac{1}{3}R_{ikmj}(p')$, $M_{km}=\frac{1}{3}R_{ij}(p')$
and $N^k_i=-\frac{2}{3}R_{ik}(p')$.

We prove the following lemma.
\begin{lem}\label{l4.1}
We have
\begin{equation}\label{4.1}
\frac{1}{18}R_{ij}(p')R_{km}(p')\alpha^{ijkm}+N_{ijk}^m\alpha_m^{ijk}+M_{ijkm}\alpha^{ijkm}=4K_{ijkm}\alpha^{ijkm}.
\end{equation}
\end{lem}

\proof We have, for any small $t
> 0$,
$$\begin{array}{lll}
     \dint_{B_t}\Delta_g r^2dV_g&=&
         \dint_{B_t}(8-\frac{2}{3}R_{ij}x^{ij}+2M_{ijk}x^{ijk}+2M_{ijkm}x^{ijkm}+2N^k_{ij}x^{ij}_k+2N_{ijk}^p
         x^{ijk}_p )\\[\mv]
     &&\times(1-\frac{1}{6}R_{ij}x^{ij}+K_{ijk}x^{ijk}+K_{ijkm}x^{ijkm})dx+ o(t^8)\\
     &=&  4\pi^2 t^4-2R_{ij}\alpha^{ij}\times2\pi^2\frac{t^6}{6}\\[\mv]
    &&+(\frac{1}{9}R_{ij}R_{km}\alpha^{ijkm}+2M_{ijkm}\alpha^{ijkm}+2N_{ijk}^p\alpha_p^{ijk}+8K_{ijkm}\alpha^{ijkm})
     2\pi^2\frac{t^8}{8} + o(t^8),
  \end{array}$$
on the other hand, we have
$$\begin{array}{lll}
   \dint_{\partial B_t}2rds_g&=&
     \dint_{\partial B_t}2r(1-\frac{1}{6}R_{ij}x^{ij}+K_{ijkm}x^{ijkm} + O(r^5) )ds_0\\[\mv]
 &=& 4 \pi^2t^4 -4\pi^2\frac{R_{ij}}{6}\alpha^{ij}t^6+2K_{ijkm}\alpha^{ijkm}2\pi^2t^8 + o(t^8).
 \end{array}$$
Now the conclusion follows from the Stokes' theorem.

$\hfill\Box$\\

Note that locally, we may write (see Lemma \ref{6.1} in the
appendix),
$$G_{p'}=-2\log r+S,$$
with
$$S=S_0(p')+a_ix^i+\frac{a_{ij}}{2}x^{ij}+O(r^{2+\alpha}).$$

We define
$$\varphi_\epsilon=
                -\log(1+\lambda|\frac{x}{\epsilon}|^2)
                        +C_\epsilon+\mu |x|^2,\s x\in B_{L\epsilon}$$
where
$$\mu=-\frac{1}{L^2\epsilon^2(1+\lambda L^2)},\s \lambda=\frac{\sqrt{3\tilde{Q}(p')}}{12}$$
and
$$C_\epsilon=\log(1+\lambda L^2)-2\log L\epsilon-\mu L^2\epsilon^2.$$

We set
$$\phi_\epsilon=\left\{
  \begin{array}{ll}
      G+\varphi_\epsilon+2\log r&x\in B_{L\epsilon}\\[\mv]
      G&x\notin B_{L\epsilon},
  \end{array}\right.$$
then, in $B_{L\epsilon}$, we have
$$\phi_\epsilon=-\log(1+\lambda
                       |\frac{x}{\epsilon}|^2)
                        +C_\epsilon+S+
                         \mu |x|^2.$$
Hence, it is easy to check that $\phi_\epsilon\in W^{2,p}(M)$ for any $p>0$.\\

We write
\begin{align*}
II(\phi_\epsilon) : & ~ = ~ \dint_M\lan \phi_\epsilon,
\phi_\epsilon\ran dV_g + 4 \int_M Q_g \phi_\epsilon dV_g - 8 \pi^2
\log \int_M \tilde{Q} e^{4
\phi_\epsilon} dV_g\\[2mm]
& = II_1 + II_2 + II_3
\end{align*}
First we will calculate the term $II_3$. In the small neighborhood
around the point $p'$, we set
$$\tilde{Q}=\tilde{Q}(p')+b_ix^i+\frac{b_{ij}}{2}x^{ij}+O(r^3),$$
then we have
$$\begin{array}{lll}
    \tilde{Q}e^{4\phi_\epsilon}\sqrt{|g|}
      &=&\frac{e^{4C_\epsilon+4S_0}}{\epsilon^4
        (1+\lambda|\frac{x}{\epsilon}|^2)^4}[(1+4a_ix^i+2a_{ij}x^{ij}+8a_ia_jx^{ij}
        +4\mu r^2)\tilde{Q}(p')+b_ix^i+\frac{b_{ij}}{2}x^{ij}+4a_ib_ix^{ij}\\[\mv]
      &&+O(r^{2+\alpha})+O(\frac{r^2\epsilon^2}{L^8})](1-\frac{R_{ij}x^{ij}}{6}
        +O(r^3))\\[\mv]
      &=&\frac{e^{4C_\epsilon+4S_0}}{\epsilon^4
        (1+\lambda|\frac{x}{\epsilon}|^2)^4}[(1+4a_ix^i+2a_{ij}x^{ij}+8a_ia_jx^{ij}
        +4\mu r^2-\frac{R_{ij}x^{ij}}{6})\tilde{Q}(p')+b_ix^i+\frac{b_{ij}}{2}x^{ij}+4a_ib_ix^{ij}\\[\mv]
      &&+O(r^{2+\alpha})+O(\frac{r^2}{L^8})].
  \end{array}$$
Therefore, by using the symmetry of the ball and the fact that $\alpha_{ij}=\frac{1}{4}\delta_{ij}$, we have
$$\begin{array}{lll}
   \dint_{B_{L\epsilon}}\tilde{Q}e^{4\phi_\epsilon}\sqrt{|g|}dV_g&=&2\pi^2
       e^{4C_\epsilon+4S_0(p')}\epsilon^4\dint_0^L\frac{1}{(1+\lambda r^2)^4}
              [\tilde{Q}(p')(1+\epsilon^2r^2(\sum\limits_i(\frac{a_{ii}}{2}+2a_i^2)+
                4\mu-\frac{R(p')}{24})\\[\mv]
   &&+\sum\limits_i(a_ib_i+\frac{b_{ii}}{8})\epsilon^2r^2+O(\epsilon r)^{2+\alpha}
   +O(\frac{r^2}{L^4})]r^3dr.
  \end{array}$$
A direct calculation then yields that
\begin{align*}
 2\pi^2\dint_0^L\frac{r^3dr}{(1+\lambda r^2)^4}
 & =\frac{\pi^2}{6\lambda^2}+O(\frac{1}{L^4}),\\[\mv]
 2\pi^2\dint_0^L\frac{r^5dr}{(1+\lambda r^2)^4} &=
\frac{\pi^2}{3\lambda^3} +O(\frac{1}{L^2}),
\end{align*}
 and
$$4\mu\epsilon^2\times2\pi^2\dint_0^L\frac{r^5dr}{(1+\lambda r^2)^4}
  =O(\frac{1}{L^4}).$$
Hence we get
$$\begin{array}{lll}
 \dint_{B_{L\epsilon}}\tilde{Q}e^{4\phi_\epsilon}\sqrt{|g|}dx
   &=&e^{4C_\epsilon+4S_0}\epsilon^4[8\pi^2
       -\frac{24\pi^2}{\lambda^2L^4}+ \frac{\pi^2}{3\lambda^3}\epsilon^2
        (\sum\limits_i(\frac{a_{ii}}{2}+2a_i^2)\tilde{Q}(p')-\frac{R(p')}{24}\tilde{Q}(p')\\[\mv]
   &&+\sum\limits_i(a_{i}b_i+\frac{b_{ii}}{8}))+O(\frac{1}{L^4})+
     O(\epsilon^{2+\alpha})+O(\frac{\epsilon^2}{L^2})].
 \end{array}$$
On the other hand, it is not difficult to check that
$$\begin{array}{lll}
   \dint_{M\setminus
   B_{L\epsilon}}\tilde{Q}e^{4\phi_\epsilon}\sqrt{|g|}dx
      &=&\dint_{L\epsilon}^{\delta}\tilde{Q}(p')\frac{e^{4S_0}}{r^5}2\pi^2dr+
         O(\frac{1}
        {L^2\epsilon^2})\\[\mv]
      &=&e^{4C_\epsilon+4S_0}\epsilon^4(\frac{24\pi^2}{\lambda^2L^4}+O(\frac{\epsilon^2}
        {L^2})).
      \end{array}$$
In sum, we have
\begin{equation}\label{4.2}
\begin{array}{lll}
     8\pi^2\log\dint_M
       \tilde{Q}e^{4\phi_\epsilon}\sqrt{|g|}dx&=&8\pi^2
       [\log8\pi^2+4(C_\epsilon
       +\log\epsilon+S_0)]\\[\mv]
     &&+\frac{\pi^2}{3\lambda^3}[\tilde{Q}(p')\sum\limits_i(\frac{a_{ii}}{2}+2a_i^2)
   +\sum\limits_i(a_{i}b_i+\frac{b_{ii}}{8})-\frac{R(p')}{24}\tilde{Q}(p')]\epsilon^2\\[\mv]
     &&+O(\epsilon^{2+\alpha})+O(\frac{\epsilon^2}{L^2})+O(\frac{1}{L^4}).
  \end{array}
\end{equation}
\\

The next, we calculate $II_1$: First of all, we have
\begin{equation}\label{4.3}
  \begin{array}{lll}
      \dint_M\lan \phi_\epsilon,\phi_\epsilon\ran dV_g&=&
        \dint_M\lan G,\phi_\epsilon\ran dV_g+\dint_{B_{L\epsilon}}\lan \varphi_\epsilon
          +2\log r, \phi_\epsilon\ran dV_g\\[\mv]
      &=&16\pi^2(C_\epsilon+S_0(p'))-2\dint_MQ\phi_\epsilon dV_g
          +\dint_{B_{L\epsilon}}\lan \varphi_\epsilon+2\log r,\varphi_\epsilon+S\ran
          dV_g.
      \end{array}
\end{equation}

We set $\eta$ to be a cut-off function which is 0 at 1 and 1 in
$[0,1/4]$ with $\eta'(1)=1$, and
$$h_\tau=\left\{
     \begin{array}{ll}
       \eta({\frac{|x|}{\tau}})+\log \tau &|x|\leq \tau \\[\mv]
       \log r&|x|\geq\tau.
     \end{array}\right.$$
Then for  fixed $\epsilon$ and $L$, we have
$$\lim_{\tau \rightarrow 0}\dint_{B_{L\epsilon}}
\lan \varphi_\epsilon+2h_\tau,\varphi_\epsilon+S\ran dV_g
=\dint_{B_{L\epsilon}}\lan \varphi_\epsilon+2\log
r,\varphi_\epsilon+S\ran dV_g.$$
On the other hand, we have
$$\begin{array}{lll}
    \dint_{B_{L\epsilon}}
       \lan \varphi_\epsilon+2h_\tau,\varphi_\epsilon+S\ran dV_g&=&
        \dint_{B_{L\epsilon}}
         \lan \varphi_\epsilon+2h_\tau,G\ran dV_g+\dint_{B_{L\epsilon}}
          \lan \varphi_\epsilon+2h_\tau,\varphi_\epsilon+2\log r\ran dV_g\\[\mv]
     &=&16\pi^2C_\epsilon+32\pi^2\eta(0)+32\pi^2\log\tau -2\dint_{B_{L\epsilon}}
       Q_g(\varphi_\epsilon+2h_\tau)\\[\mv]
     &&+\dint_{B_{L\epsilon}}\lan\varphi_\epsilon,\varphi_\epsilon\ran dV_g
     +\dint_{B_{L\epsilon}}\lan\varphi_\epsilon,2\log r+2h_\tau \ran dV_g\\[\mv]
     &&+\dint_{B_{L\epsilon}}\lan 2\log r,2h_\tau \ran dV_g.
\end{array}$$

Therefore we get
\begin{equation}\label{4.4}
\begin{array}{l}
   \dint_{B_{L\epsilon}}\lan \varphi_\epsilon+2\log r,\varphi_\epsilon+S\ran dV_g\\[\mv]
 \begin{array}{lll}
    &=&
     32\pi^2\eta(0)-2\dint_{B_{L\epsilon}}Q_g(\varphi_\epsilon+2\log r)
     +\dint_{B_{L\epsilon}}\lan\varphi_\epsilon,\varphi_\epsilon\ran dV_g\\[\mv]
     &&+\dint_{B_{L\epsilon}}\lan\varphi_\epsilon,4\log r\ran dV_g
     +\lim\limits_{\tau\rightarrow 0}(\dint_{B_{L\epsilon}}\lan 2\log r,2h_\tau\ran dV_g
     +32\pi^2\log\tau)\\[\mv]
    &=&
     32\pi^2\eta(0)-2\dint_{B_{L\epsilon}}Q_g(\varphi_\epsilon+2\log r)
     +\dint_{B_{L\epsilon}}\Delta_g\varphi_\epsilon\Delta_g\varphi_\epsilon dV_g\\[\mv]
     &&+4\dint_{B_{L\epsilon}}\Delta_g\varphi_\epsilon\Delta_g\log r dV_g
     +\lim\limits_{\tau\rightarrow 0}(\dint_{B_{L\epsilon}}\Delta_g 2\log r
       \Delta_g2h_\tau dV_g
     +32\pi^2\log\delta)\\[\mv]
     &&+\dint_{B_{L\epsilon}}\frac{2}{3}R\lan d(\varphi_\epsilon+2\log r),
      d(\varphi_\epsilon+2\log r)\ran dV_g\\[\mv]
     && -\dint_{B_{L\epsilon}}2Ric(d(\varphi_\epsilon+2\log r),
      d(\varphi_\epsilon+2\log r))dV_g.
  \end{array}
\end{array}
\end{equation}

By a simple calculation, one gets
\begin{equation}\label{4.5}
\begin{array}{lll}
   \dint_{B_\tau}(\Delta_g 2\log r) \Delta_g (2h_\tau) d V_g&=&
          \dint_{B_\tau}\Delta_0 (2\log r) \Delta_0 (2\eta(\frac{|x|}{\tau}))dx+ O(\tau)\\[\mv]
    &=&-32\pi^2\eta(0)+16\pi^2+ O(\tau).
  \end{array}
\end{equation}

To compute $\dint_{B_{L\epsilon}\setminus B_\delta}\Delta_g\log
r\Delta_g\log r$, we first verify that, for any smooth function
$f$, $g$ which are smooth in $(t_0, t_1)$, we have
$$\begin{array}{lll}
     \Delta_gf(r)&=&(\delta_{km}+M_{ij}^{km}x^{ij}+M_{ijs}^{km}x^{ijs}+M_{ijst}^{km}x^{ijst}+O(r^5))
        (f''\frac{x_{km}}{r^2}+f'\frac{\delta_{km}}{r}-f'\frac{x_{km}}{r^3})
        +N^k\frac{x_k}{r}f'\\[\mv]
     &=&f''+f'(\frac{3}{r}-\frac{R_{ij}x^{ij}}{3r}+
         \frac{M_{ijk}x^{ijk}+N_{ij}^kx_k^{ij}}{r}+\frac{M_{ijkm}x^{ijkm}+N_{ijk}^mx_m^{ijk}}{r})+O(r^5|f''|)+O(r^4|f'|).
  \end{array}$$
Here, we use the fact that
$M_{ij}^{km}x_{km}^{ij}=M_{ijst}^{km}x_{km}^{ijst}=0$. Then,
applying Lemma 4.1, we get
\begin{equation}\label{4.6}
\begin{array}{l}
\dint_{B_{t_1}\setminus B_{t_0}}\Delta_g f(|x|)\Delta_g g(|x|)dV_g\\[\mv]
  \begin{array}{lll}
    &=&\dint_{t_0}^{t_1}f''g''(1-\frac{R}{24}r^2+K_{ijkm}\alpha^{ijkm}r^4)2\pi^2r^3dr\\[\mv]
    &&+\dint_{t_0}^{t_1}(f'g''+f''g')\frac{1}{r}(3-\frac{5R}{24}r^2+7K_{ijkm}\alpha^{ijkm}r^4)2\pi^2r^3dr\\[\mv]
    &&+\dint_{t_0}^{t_1}f'g'\frac{1}{r^2}(9+33K_{ijkm}\alpha^{ijkm}r^4-\frac{7R}{8}r^2+\frac{1}{9}R_{ij}R_{km}
       \alpha^{ijkm}r^2)2\pi^2r^3dr\\[\mv]
    &&+\dint_{t_0}^{t_1}\left(O(r^8|f''g''|)+O(r^7(|f''g'|+|f'||g''|))+O(r^6|f'g'|)\right)\\[\mv]
    &=&\dint_{t_0}^{t_1}(f''g''+(f'g''+f''g')\frac{3}{r}+f'g'\frac{9}{r^2})2\pi^2r^3\\[\mv]
    &&+R\dint_{t_0}^{t_1}(-f''g''\frac{r^2}{24}-\frac{5r}{24}(f'g''+f''g')-\frac{7}{8}f'g')2\pi^2r^3\\[\mv]
    &&+K_{ijkm}\alpha^{ijkm}\dint_{t_0}^{t_1}(f''g''r^4+7(f'g''+f''g')r^3+33f'g'r^2)2\pi^2r^3dr\\[\mv]
    &&+R_{ij}R_{km}\alpha^{ijkm}\dint_{t_0}^{t_1}\frac{1}{9}f'g'r^2
      2\pi^2r^3dr\\[\mv]
     &&+\dint_{t_0}^{t_1}\left(O(r^8|f''g''|)+O(r^7(|f''g'|+|f'||g''|))+O(r^6|f'g'|)\right)dr.
  \end{array}
\end{array}
\end{equation}

Then, choosing $f=g=2\log r$, $t_1=L\epsilon$, $t_0=\tau$, we get
\begin{equation}\label{4.7}
\begin{array}{lll}
   \dint_{B_{L\epsilon}\setminus B_\tau}\Delta_g (2\log r) \Delta_g (2h_\tau) dV_g&=&
          \dint_{B_{L\epsilon}\setminus B_\tau}\Delta_g (2\log r) \Delta_g (2\log r) dV_g\\[\mv]
    &=& 40K_{ijkm}\alpha^{ijkm}
       \pi^2(L\epsilon)^4+\frac{2\pi^2}{9}R_{ij}R_{km}\alpha^{ijkm}(L\epsilon)^4\\[\mv]
    &&-2R\pi^2(L\epsilon)^2
      +32\pi^2\log L\epsilon-32\pi^2\log\tau\\[\mv]
   && + O(\tau)+O(L\epsilon)^5.
  \end{array}
\end{equation}
\\

Now we will calculate the term $\int_{B_{L\epsilon}}\Delta_g
\varphi_\epsilon \Delta_g(\varphi_\epsilon+4\log r) dV_g$: In
(\ref{4.6}), we choose $f=\varphi_\epsilon$,
$g=\varphi_\epsilon+4\log r$, $t_0=0$, $t_1=L\epsilon$ then we get
\begin{equation}\label{4.8}
\begin{array}{lll}
   \dint_{B_{L\epsilon}}\Delta_g
          \varphi_\epsilon \Delta_g(\varphi_\epsilon+4\log r) dV_g&=&
            -\frac{88}{3}\pi^2+\frac{16\pi^2}{\lambda L^2}-16\pi^2\log(1+\lambda L^2)\\[\mv]
            &&-R\epsilon^2\frac{8\pi^2}{9\lambda}+2\pi^2R(L\epsilon)^2\\[\mv]
            &&-40K_{ijkm}\alpha^{ijkm}\pi^2(L\epsilon)^4-\frac{2\pi^2}{9}R_{ij}R_{km}
            \alpha^{ijkm}(L\epsilon)^4\\[\mv]
           &&+O(\epsilon^4L^2)+\frac{\epsilon^2}{L^2}+O(L\epsilon)^5.
\end{array}
\end{equation}

By a direct calculation, we have
\begin{equation}\label{4.9}
\begin{array}{l}
   \dint_{B_{L\epsilon}}\frac{2}{3}R ( \nabla_g(\varphi_\epsilon+2\log r),
         \nabla_g(\varphi_\epsilon+2\log r))dV_g\\[\mv]
   \s\s\s\s=\frac{2}{3}\dint_0^{L\epsilon}R(p')(\frac{2\epsilon^2}{(\epsilon^2+\lambda r^2)r}+2\mu r)^22\pi^2r^3\\[\mv]
    \s\s\s\s\s+\frac{2}{3}
         \dint_{B_{L\epsilon}}(R_{,i}(p')x^i+O(r^2)
         )(\frac{2\epsilon^2}{(\epsilon^2+\lambda r^2)r}+2\mu r)^2
         (1+O(r^3))dx\\[\mv]
   \s\s\s\s=\frac{8}{3\lambda}R(p')\pi^2\epsilon^2+
         \dint_{B_{L\epsilon}}(\frac{2\epsilon^2}{(\epsilon^2+\lambda r^2)r}+2\mu r)^2
         O(r^2)dx\\[\mv]
   \s\s\s\s=\frac{8}{3\lambda}R(p')\pi^2\epsilon^2+O(\epsilon^4L^2)+O(\frac{\epsilon^2}{L^2}),
  \end{array}
\end{equation}
and
\begin{equation}\label{4.10}
\begin{array}{l}
   \dint_{B_{L\epsilon}}2Ric(\nabla_g(\varphi_\epsilon+2\log r),
                  \nabla_g(\varphi_\epsilon+2\log r))]dV_g\\[\mv]
   \s\s\s\s=\frac{1}{2}R(p')\dint_0^{L\epsilon}(\frac{2\epsilon^2}{(\epsilon^2+\lambda r^2)r}+2\mu r)^22\pi^2r^3dr\\[\mv]
   \s\s\s\s\s\s+2
         \dint_{B_{L\epsilon}}g^{is}g^{jt}(R_{ij,k}(p')x^k
         +O(r^2))(\frac{2\epsilon^2}{(\epsilon^2+\lambda r^2)r^2}+2\mu )^2x_{st}
         (1+O(r^3))dx\\[\mv]
   \s\s\s\s=\frac{2}{\lambda}R(p')\pi^2\epsilon^2+2
         \dint_{B_{L\epsilon}}(R_{ij,k}(p')x^k
         +O(r^2))(\frac{2\epsilon^2}{(\epsilon^2+\lambda r^2)r^2}+2\mu )^2x^{ij}
         (1+O(r^3))dx\\[\mv]
   \s\s\s\s=\frac{2}{\lambda}R(p')\pi^2\epsilon^2+
         \dint_{B_{L\epsilon}}(\frac{2\epsilon^2}{(\epsilon^2+\lambda r^2)r^2}+2\mu)^2
         O(r^4)dx\\[\mv]
   \s\s\s\s=\frac{2}{\lambda}R(p')\pi^2\epsilon^2+O(\epsilon^4L^2)+O(\frac{\epsilon^2}{L^2}).
  \end{array}
\end{equation}
Together with (\ref{4.3})-(\ref{4.5}) and
(\ref{4.7})-(\ref{4.10}),  we obtain the following identity
\begin{equation}
\begin{array}{lll}
     II_\epsilon(u_\epsilon)&=&  II_1 + II_2 + II_3\\[\mv]
     &=& -16\pi^2\log\lambda-8\pi^2\log8\pi^2
      +\frac{8\pi^2}{3}-16\pi^2+2\dint_MQG-16\pi^2S_0\\[\mv]
  &&-\frac{\epsilon^2\pi^2}{3\lambda^3}
    (\tilde{Q}(p')\sum\limits_{i}(\frac{a_{ii}}{2}+2a_i^2)
     +\sum_{i}(a_ib_i+\frac{b_{ii}}{8})-\frac{R(p')}{36}\tilde{Q}(p'))\\[\mv]
  &&+O(\frac{\epsilon^2}{L^2})+O(\epsilon^{2+\alpha})+O(\frac{1}{L^4})+O(\epsilon^4L^2)
  +O((L\epsilon)^5).
  \end{array}
\end{equation}

Proof of Theorem \ref{main2} : we set $L =
\frac{\log \frac{1}{\epsilon}}{\epsilon^{\frac{1}{2}}}$, then
$$\epsilon^2  \gg O(\frac{\epsilon^2}{L^2})+O(\epsilon^{2+\alpha})+O(\frac{1}{L^4})+O(\epsilon^4L^2)
  +O((L\epsilon)^5)$$
when $\epsilon$ is very small.
Therefore, we get Theorem \ref{main2}. $\hfill\Box$

\section{The conformal case}
In this section, we will discuss the local conformal flat case of
Theorem \ref{main2}.

In this situation, locally we may write
$$g=e^{2f}\sum_idx^i\otimes dx^i\s with \s f=c_ix^i+\frac{1}{2}c_{ij}x^{ij}+O(r^3),$$
and
$$\tilde{Q}=\tilde{Q}(p')+b_ix^i+\frac{1}{2}b_{ij}x^{ij}+O(r^3).$$
Note that by the conformal property of $P_g$, the corresponding
Green function have the following local expression:
$$G=-2\log|x|+S_0(p')+a_ix^i+\frac{1}{2}a_{ij}x^{ij}+O(r^3).$$

When $f=0$, we can use Theorem \ref{main2} to obtain: if
$$\sum_i(\frac{a_{ii}}{2}+2a_i^2+\frac{1}{\tilde{Q}(p')}(
a_ib_i+\frac{b_{ii}}{8}))>0,$$ then (\ref{003}) has a solution.

For the general case, we set $g'=e^{-2f}g$, then applying Lemma
\ref{l3.5}, we get $G'_{p'}=G+f$, and then
$$a_i'=a_i+c_i,\s and \s a_{ii}'=a_{ii}+c_{ii}.$$\\
Thus we have the following results

\begin{thm}\label{t5.1}
Let $(M,g)$ be a close 4-dimensional manifold with $k = 8\pi^2$
and $P_g$ is positive. Suppose further that it is locally
conformal flat near $p'$. If
$$ \sum_i \frac{a_{ii}+c_{ii}}{2}+2(a_i+c_i)^2+\frac{1}{\tilde{Q}(p')}(
(a_i+c_i)b_i+\frac{b_{ii}}{8})>0,$$ then
equation (\ref{003}) has a minimal solution.\\
\end{thm}

As a corollary, we have\\

\begin{cor}\label{c5.2}
With the same assumption as in Theorem 5.1. If
$$\sum_i \frac{a_{ii}+c_{ii}}{2}+2(a_i+c_i)^2>0,$$
then in the conformal class of $(M,g)$ there is a constant $Q$-curvature.\\
\end{cor}

To end this section, we propose the following conjecture:\\

\noindent{\bf Conjecture:} {\it  Let $(M,g)$ be a locally
conformal flat closed Riemannian manifold of dimension four, with
$k = 8\pi^2$ and $P_g$ is positive. Then we have
$$\sum_{i}(\frac{a_{ii}+c_{ii}}{2}
     +2(a_i+c_i)^2) \geq 0, ~~ \mbox{at the point}
     ~ p'~{\rm where}~ \Lambda_g(p') = \min_{x\in M} \Lambda_g(8\pi^2,x),$$
and the equality holds if and only if  $(M,g)$
is in the conformal class of the standard 4-sphere.}\\

Let $\tilde{g}=e^{2G}g$, then we have
$$Q_{\tilde{g}}(x)=0$$
for any $x\neq p$. Near p, we can write
$$\tilde{g}=\frac{e^{S_0(p)+(c_i+a_i)x^i+(c_{ij}+a_{ij})x^{ij}}}{r^2}=
\frac{e^{S_0(p)}}{r^2}(\theta_ix^i+\theta_{ij}x^{ij}+O(|x|^3)).$$
So the above conjecture is equivalent to that
$$\sum_i\theta_{ii}> 0$$
when $M\neq S^4$. So, this problem is very similar to the positive
mass problem.

\section{Appendix}
{\small Suppose $Ker{P_g}=\{constant\}$.
Let $G$ be the Green function which satisfies
$$P_gG+2Q_g=16\pi^2\delta_p.$$
As a corollary of a result in \cite{N}, we have the following

\begin{lem}\label{6.1}
In a normal coordinate system of $p$, we
have
$$G=-2\log r+S_0+a_ix^i+a_{ij}x^{ij}+O(r^{2+\alpha}).$$
\end{lem}

However, for the reader's sake, we give a brief proof of this Lemma here:

\proof In a normal coordinate system, we set
$$|g|=1-\frac{1}{3}R_{ij}x^{ij}+O(r^3),\s and\s
g^{km}=\delta^{km}-
\frac{1}{3}R_{kijm}x^{ij}+O(r^3)$$ where
$\varphi_{ijk}$ and $\theta_{ijk}$ are smooth.

Given a smooth function $F$, we have
$$\begin{array}{lll}
  \Delta_gF(|x|)&=&
 \frac{1}{\sqrt{|g|}}\frac{\partial}{\partial x^k}
  (\sqrt{|g|}g^{km}\frac{\partial }{\partial x^m}F)\\[\mv]
  &=&\frac{\partial}{\partial x_k}(g^{km}F'\frac{x_m}{r})+\frac{1}{2}g^{km}F_m\frac{\partial}{\partial x_k}
  \log|g|\\[\mv]
  &=&\frac{\partial}{\partial x_k}(F'\frac{x_k}{r}-\frac{1}{3}R_{kijm}F'\frac{x^{kij}}{r}+F'O(r^3))
     -\frac{1}{3}R_{ij}F'\frac{x^{ij}}{r}+O(F'r^2)\\[\mv]
  &=&\frac{\partial}{\partial x_k}(F'\frac{x_k}{r}+F'O(r^3))
     -\frac{1}{3}R_{ij}F'\frac{x^{ij}}{r}+O(F'r^2)\\[\mv]
    &=&\Delta_0F
     -\frac{1}{3}R_{ij}F'\frac{x^{ij}}{r}+O(F'r^2)+O(F''r^3).
 \end{array}$$

Then
$$\Delta_g(-2\log r)=-\frac{4}{r^2}+\frac{2}{3}R_{ij}\frac{x^{ij}}{r^2}+O(r)$$
and
$$\Delta_g(-\frac{4}{r^2})=\Delta_0(-\frac{4}{r^2})-\frac{8R_{ij}x^{ij}}{3r^4}
   +O(\frac{1}{r})=16\pi^2\delta_0-\frac{8R_{ij}x^{ij}}{3r^4}+O(\frac{1}{r}).$$
It is easy to check that
$$\Delta_g\frac{2}{3}R_{ij}\frac{x^{ij}}{r^2}=\Delta_0\frac{2}{3}R_{ij}\frac{x^{ij}}{r^2}
 +O(\frac{1}{r})
=\frac{4R}{3r^2}-\frac{16R_{ij}x^{ij}}{3r^4}.$$
Hence, we get
$$\Delta_g^2(-2\log r)=16\pi^2\delta_p+\frac{4R}{3r^2}-8\frac{R_{ij}x^{ij}}{r^4}
+O(\frac{1}{r}).$$

Moreover, we have
$$\begin{array}{lll}
div(\frac{2}{3}R_g(-d2\log r)-2Ric_g\lan
d(-2\log r),\cdot\ran)&=&
      \frac{2}{3}R_p(p')(2\log r)_{kk}-2R_{km}(p')(2\log r)_{km}+O(\frac{1}{r})\\[\mv]
 &=&\frac{2}{3}R_g(p')\frac{4}{r^2}-4R_g(p')\frac{1}{r^2}+8R_{km}\frac{x^{km}}{r^4}+O(\frac{1}{r}).
\end{array}$$
We therefore have
$$P_g(-2\log r)=16\pi^2\delta_0+O(\frac{1}{r}).$$

We set
$$G=-2\log r+S$$
where $S\in C^{1,\alpha}$. Then, we get
$$\Delta_g^2S=P_gS+O(\frac{1}{r})=P_gG+2P_g\log r+O(\frac{1}{r})=O(\frac{1}{r}).$$
This proves the lemma. $\hfill\Box$

\noindent Jiayu Li\\
{ICTP, Mathematics Section, Strada Costiera 11, 34014 Trieste,
Italy, and Academy of Mathematics and Systems Sciences, Chinese
Academy of Sciences, Beijing 100080, P.R. China}\\
{\it E-mail address:} jyli@ictp.it\\

\noindent Yuxiang  Li\\
{ICTP, Mathematics Section, Strada Costiera 11, 34014 Trieste,
Italy}\\
{\it E-mail address:} liy@ictp.it\\

\noindent Pan Liu\\
{Department of Mathematics, East China Normal University, 3663,
Zhong Shan North Rd,
Shanghai 200062, P.R. China\\
{\it E-mail address:} pliu@math.ecnu.edu.cn\\

\end{document}